\def\docdate{May 12, 2006}

\documentclass{ltxguide}

\title{Closure properties in positively decreasing and related\\ distributions under dependence}

\author{Dimitrios G. Konstantinides,\,Charalampos  D. Passalidis,\,Nikolaos E. Porichis} 

\date{\docdate}

\DeclareMathAccent{\wtilde}{\mathord}{largesymbols}{"65}

\raggedright
\advance\textheight24pt
\usepackage{amsmath}
\usepackage{amssymb}
\usepackage{amsthm}
\usepackage{latexsym}
\usepackage{amssymb}

\newtheorem{theorem}{Theorem}[section]
\newtheorem{lemma}{Lemma}[section]
\newtheorem{corollary}{Corollary}[section]
\newtheorem{remark}{Remark}[section]
\newtheorem{definition}{Definition}[section]
\newtheorem{proposition}{Proposition}[section]
\newtheorem{example}{Example}[section]
\newtheorem{assumption}{Assumption}[section]
\numberwithin{equation}{section}
\newcommand{\bth}{\begin{theorem}}
\newcommand{\ethe}{\end{theorem}}

\newcommand{\bre}{\begin{remark}}
\newcommand{\ere}{\end{remark}}

\newcommand{\ble}{\begin{lemma}}
\newcommand{\ele}{\end{lemma}}
\newcommand{\bde}{\begin{definition}}
\newcommand{\ede}{\end{definition}}

\newcommand{\bco}{\begin{corollary}}
\newcommand{\eco}{\end{corollary}}

\newcommand{\bpr}{\begin{proposition}}
\newcommand{\epr}{\end{proposition}}

\newcommand{\bexer}{\begin{exercise}}
\newcommand{\eexer}{\end{exercise}}
\newcommand{\breh}{\begin{hint}}
\newcommand{\ereh}{\end{hint}}

\newcommand{\halmos}{\hfill \qed}

\newcommand{\bexam}{\begin{example}}
\newcommand{\eexam}{\end{example}}

\newcommand{\pr} {{\bf Proof.}}

\newcommand{\bfi}{\begin{fig}}
\newcommand{\efi}{\end{fig}}

\newcommand{\beao}{\begin{eqnarray*}}
\newcommand{\eeao}{\end{eqnarray*}\noindent}

\newcommand{\beam}{\begin{eqnarray}}
\newcommand{\eeam}{\end{eqnarray}\noindent}
\newcommand{\E}{\mathbf{E}}
\newcommand{\PP}{\mathbf{P}}

\newcommand{\xto}{x\to\infty}

\newcommand{\bH}{\overline{H}}
\newcommand{\bF}{\overline{F}}

\newcommand{\bG}{\overline{G}}

\newcommand{\bbr}{{\mathbb R}}

\newcommand{\bbn}{{\mathbb N}}

\newcommand{\vep}{\varepsilon}

\begin{document}
\date{{\small \today}}

\maketitle
\begin{abstract}
We consider closure properties in the class of positively decreasing distributions. Our results stem from different types of dependence, but each type belongs in the family of asymptotically independent dependence structure. Namely we examine the closure property with respect to minimum, maximum, convolution product and convolution. Furthermore, we take into account some closure properties of the class of generalized subexponential positively decreasing distributions, as also we introduce and study the class of the generalized long-tailed positively decreasing distributions. Also we revisited the (independent) convolution closure problem of subexponentiality, in the case of subexponential positively decreasing class. In some classes we discuss the closedness of randomly stopped sums. In the last section we study the closure property with respect to minimum for two classes of random vectors.
\end{abstract}
\vspace{3mm}

\section{Introduction} \label{sec.KP.1}

Let consider two non-negative random variables $X$ and $Y$, following distributions $F$ and $G$ respectively, with supports the positive half-axis $\bbr_+$. We denote by $\bF:=1-F$ the distribution tail, hence $\bF(x)=\PP[X>x]$ and holds $\bF(x)>0$
for any $x \geq 0$, except if it is said differently. For two positive functions $f(x)$ and $g(x)$, the asymptotic relation $f(x)=o[g(x)]$, as $\xto$, means
\beao
\limsup_{\xto}\dfrac{f(x)}{g(x)} = 0\,,
\eeao
the asymptotic relation $f(x)=O[g(x)]$, as $\xto$, holds if
\beao
\limsup_{\xto} \dfrac{f(x)}{g(x)} < \infty\,.
\eeao
while the asymptotic relation $f(x)\asymp g(x)$, as $\xto$, holds if  both $f(x)=O[g(x)]$ and $g(x)=O[f(x)]$ are true. For a real number $x$, we denote $x^{+}:=\max\{x,0\}$ and for two real numbers $x$ and $y$ we write $x\vee y :=\max\{x,\,y\}$ and $x\wedge y:=\min\{x,\,y\}$.

Now, let introduce some classes of distribution with heavy tails. We say that distribution $F$ is heavy-tailed, symbolically $F \in \mathcal{H}$ if
\beao
\int_{-\infty}^{\infty} e^{\vep \,y} \,F(dy) =\infty\,,
\eeao 
for any $\vep >0$. A very large subclass of this class of heavy-tailed distributions is the class of long tailed distributions, symbolically $F\in \mathcal{L}$, whence for any (or equivalently for some) $y\neq 0$ holds
\beao
\lim_{\xto} \dfrac{\bF(x-y)}{\bF(x)}=1\,.
\eeao

Another well-known class, of heavy-tailed distributions, is the class of subexponential distributions, symbolically $F\in \mathcal{S}$, which is characterized by the property
\beao
\lim_{\xto}\dfrac{\overline{F^{*n}}(x)}{\bF(x)} = n\,,
\eeao
for any integer $n\geq 2$, where by $F^{*n}$ we denote the $n$-th order convolution power of $F$. The subexponential class was introduced in \cite{chistyakov:1964} in relation with some application on branching processes. Furthermore, it was used in many application in actuarial and financial mathematics, see for example \cite{yang:li:2019}, \cite{lu:yuan:2022}, \cite{cheng:cheng:2018} and \cite{geng:ji:wang:2019} among many others.

The following class was introduced in \cite{feller:1969} and is called class of dominatedly varying distribution, symbolically $F \in \mathcal{D}$. We say  $F \in \mathcal{D}$ if for some (or equivalently for any) $b \in (0,\,1)$ holds
\beam \label{eq.KPP.1.0}
\limsup_{\xto} \dfrac{\bF(b\,x)}{\bF(x)} < \infty\,.
\eeam
It is well known that $\mathcal{D} \nsubseteq \mathcal{S}$, $\mathcal{S} \nsubseteq \mathcal{D}$ and $\mathcal{D} \cap \mathcal{S}=\mathcal{D} \cap \mathcal{L} \neq \emptyset$ (see \cite{goldie:1978} for a proof).

The following important condition, can be found in \cite{bingham:goldie:teugels:1987}, however it was initially suggested by de Haan - Resnick (1983). We say that a distribution $F$ belongs to the class of positively decreasing distributions, symbolically $F \in \mathcal{P_D}$, if for some (or equivalently for any) $v>1$ holds
\beam \label{eq.KPP.1.1}
\limsup_{\xto} \dfrac{\bF(v\,x)}{\bF(x)} < 1\,.
\eeam 
The classes $\mathcal{D}$ and $\mathcal{P_D}$ in some sense present a symmetry. Indeed, let see the lower and upper Matuszewska indexes as follows
\beao
\beta_F:=\sup\left\{-\dfrac{\ln \bF^*(v)}{\ln v}\;:\;v>1 \right\}\,,\qquad \alpha_F:=\inf\left\{-\dfrac{\ln \bF_*(v)}{\ln v}\;:\;v>1 \right\}\,,
\eeao
where
\beao
\bF_*(v):=\liminf_{\xto} \dfrac{\bF(v\,x)}{\bF(x)}\,,\qquad \bF^*(v):=\limsup_{\xto} \dfrac{\bF(v\,x)}{\bF(x)}\,,
\eeao
for any $v>1$, and these indexes $\alpha_F\geq \beta_F$ were introduced in \cite{matuszewska:1964} and can be defined only for distributions with infinite right end point, namely when $r_F:=\sup\{y\;:\;\bF>0\}=\infty$. In general hold the inequalities $0\leq \beta_F \leq \alpha_F \leq \infty$. In relation with symmetry we can see that $F \in \mathcal{D}$ if and only if $\alpha_F < \infty$, while $F\in \mathcal{P_D}$ if and only if $\beta_F >0$. From \cite[Proposition 2.2.1]{bingham:goldie:teugels:1987} we have that for any $q\in (0,\,\beta_F)$ there exist positive constants $C$ and $x_0$ such that for any $v>1$
\beam \label{eq.KPP.1.2}
\dfrac{\bF(v\,x)}{\bF(x)} \leq C\,v^{-q}\,,
\eeam
uniformly for any $x> x_0$. Relation \eqref{eq.KPP.1.2} characterizes class $\mathcal{P_D}$ and the reason is that $F\in \mathcal{P_D}$ if and only if $\beta_F >0$, thence this relation is true for distributions from $\mathcal{P_D}$.

The following classes were introduced in \cite{chover:ney:waingr:1973(a)}, \cite{chover:ney:waingr:1973(b)}. A distribution $F$ with support on the whole real axis, belongs to the exponential type class, symbolically $F \in \mathcal{L}(\gamma)$, $\gamma\geq 0$, if for any $y\geq 0$ the following limit relation holds
\beao
\lim_{\xto} \dfrac{\bF(x-y)}{\bF(x)}= e^{ \gamma\,y}\,.
\eeao
A distribution $F$ with support $[0,\,\infty)$ or the whole real axis, belongs to the convolution equivalent class, symbolically $F \in \mathcal{S}(\gamma)\,,\;\gamma\geq 0$, if $F \in \mathcal{L}(\gamma)$ and the following limit relation holds
\beao 
\lim_{\xto} \dfrac{\overline{F^{2*}}(x)}{\bF(x)}= 2\,\widehat{F}( \gamma) < \infty\,,
\eeao
where
\beao 
\widehat{F}( \gamma) :=\int_{-\infty}^\infty e^{\gamma\,y}\,F(dy)\,.
\eeao
 
We can notice that $\mathcal{S}(\gamma) \subsetneq \mathcal{L}(\gamma)$ for any $\gamma \geq 0$. If $\gamma=0$, then $\mathcal{L}(0)=\mathcal{L}$ and $\mathcal{S}(0)=\mathcal{S}$. If $\gamma>0$ then these two classes belongs to light-tailed distributions.

Next, in \cite{konstantinides:tang:tsitsiashvili:2002} was introduced the class $\mathcal{A}:=\mathcal{S}\cap \mathcal{P_D}$, which contained the most used subexponential distributions. For the study of $\mathcal{A}$ and its application in risk theory see in \cite{chen:2011}, \cite{konstantinides:2008} and \cite{tang:yuan:2014} among others. Recently in \cite{konstantinides:passalidis:2024} was introduced the class $\mathcal{T}:=\mathcal{L}\cap \mathcal{P_D}$ and it was studied with respect to closure properties. Another class represent the intersection $\mathcal{D}\cap \mathcal{P_D}$, found in \cite{dindiene:leipus:2016} or in \cite{yang:leipus:siaulys:2016}. It is well-known that $\mathcal{D}\cap \mathcal{L}\subsetneq \mathcal{S}\subsetneq \mathcal{L}\subsetneq \mathcal{H}$.

The following classes are generalizations of the long tailed and subexponential distributions, introduced by \cite{shimura:watanabe:2005} and \cite{klueppelberg:1990} respectively. We say that distribution $F$ belongs to the class of generalized long tailed distributions, symbolically $F \in \mathcal{OL}$, if for some (or equivalently of all) $t \neq 0$ holds
\beam \label{eq.KPP.1.3}
\limsup_{\xto} \dfrac{\bF(x-t)}{\bF(x)} < \infty\,,
\eeam
and we say that distribution $F$ belongs to the class of generalized subexponential distributions, symbolically $F \in \mathcal{OS}$, if holds
\beam \label{eq.KPP.1.4}
\limsup_{\xto} \dfrac{\overline{F^{*2}}(x)}{\bF(x)} < \infty\,.
\eeam
It is remarkable that these two classes contain both heavy-tailed and light-tailed distributions. Even more, we can see that $\mathcal{OS} \subsetneq \mathcal{OL}$, $\mathcal{S} \subsetneq \mathcal{OS}$ and $\mathcal{L} \subsetneq \mathcal{OL}$.

In \cite{konstantinides:leipus:siaulys:2022} was introduced the class $\mathcal{OA}:=\mathcal{OS}\cap \mathcal{P_D}$ and was studied its closure properties with respect to convolution product. 

\bde \label{def.KPP.1.1}
We say that a distribution $F$ belongs to the class of generalized long-tailed positively decreasing distributions if $F \in \mathcal{OL}\cap \mathcal{P_D}$, symbolically $F \in \mathcal{OT}$. Namely $\mathcal{OT}:=\mathcal{OL}\cap \mathcal{P_D}$.
\ede

We can check that $\mathcal{OA}\subseteq \mathcal{OT}$ and from the inclusion $F \in \mathcal{OT}$ follows $\beta_F >0$. This new class will be examined to its closure properties. In section 2 we study the closure property of classes $\mathcal{P_D}$, $\mathcal{OL}$ and $\mathcal{OT}$ with respect to maximum and minimum, under dependence structures that include as special case the independence. Further we show the closure property of $\mathcal{OT}$ with respect to the finite mixture, with arbitrarily dependent random variables, and provide some closure properties of $\mathcal{P_D}$, $\mathcal{OA}$ and $\mathcal{OT}$ with respect to convolution of independent random variables. In section 3 we study the nesessary and sufficient condition of the closure under (independent) convolution of class $\mathcal{A}$. In section 4 we consider the class $\mathcal{P_D}$ with respect to product convolution, and based on this result we study some related classes in the same problem. In section 5 we study the closure property with respect to minimum for the classes of positenely decreasing and dominatedly varying random vectors.

\section{Closure properties with respect to convolution, mixture, maximum and minimum} \label{sec.KP.1.2}

Let consider the distribution tail of the maximum and minimum of $X$ and $Y$.
Namely we find
\beam \label{eq.KPP.2.5}
\bF_{X\vee Y} (x)= \PP[X\vee Y >x] = \bF(x) + \bG(x) - \PP[X>x\,,\;Y>x]\,,
\eeam
for the maximum and 
\beam \label{eq.KPP.2.6}
\bF_{X\wedge Y} (x)= \PP[X\wedge Y >x] = \PP[X>x\,,\;Y>x]\,,
\eeam
for the minimum, and the distribution tail of convolution of two independent random variables takes the following form
\beam \label{eq.KPP.2.7}
\overline{F*G} (x)= \PP[X + Y >x] = \int_{-\infty}^{\infty} \bF(x-y)\,G(dy)\,.
\eeam

Next, for some distribution class $\mathcal{B}$, we say that $\mathcal{B}$ is closed with respect to finite mixture if from $F,\,G \in \mathcal{B}$ follows that for any $p\in (0,\,1)$ follows $(p\,F +(1-p)\,G) \in \mathcal{B}$.

The dependence among random variables in general can not be ignored in actuarial and financial applications. For example, if $X$ and $Y$ represent claim sizes of two different portfolios in a time interval, then in most cases there is essential dependence structure. For example the motor insurance and immobility insurance can be both relevant to a severe earthquake. For dependence among main claims we refer to \cite{geluk:tang:2009}, \cite{chen:yuen:2009}, \cite{wang:2011}, \cite{geng:liu:wang:2023} and\cite{cheng:cheng:2018} for risk theory applications and to \cite{li:2022a}, \cite{li:2023} and \cite{chen:lin:2022} for risk management applications. 

In first step we consider the closure property with respect to maximum under a dependence structure, introduced in \cite{ko:tang:2008}, defined as follows.

\begin{assumption} \label{ass.KPP.1}
For $(i,\,j)\in \{(1,\,2)\,,\;(2,\,1)\}$ the asymptotic relation
\beam \label{eq.KPP.1}
\dfrac{\PP[X_j > x-t\;\mid\;X_i=t]}{\PP[X_j > x-t]} = O(1)\,,
\eeam
as $\xto$, holds uniformly for any $t\in [x_0,\,x]$ for some large enough $x_0>0$. Here, the uniformity is understood in the sense
\beam \label{eq.KPP.2}
\limsup_{\xto} \sup_{x_0 \leq t \leq x} \dfrac{\PP[X_j > x-t\;\mid\;X_i=t]}{\PP[X_j > x-t]} < \infty\,.
\eeam
\end{assumption} 

The following result establishes the closure property with respect to maximum in the classes $\mathcal{P_D}$, $\mathcal{OL}$ and $\mathcal{OT}$, under the Assumption  \ref{ass.KPP.1}.

\bth \label{th.KPP.1}
Let $X$ and $Y$ be non-negative random variables with distributions $F$ and $G$ respectively, whose supports are $\bbr_+=[0,\,\infty)$. If $X,\,Y$ satisfy Assumption \ref{ass.KPP.1}, then
\begin{enumerate}
\item
If $F \in \mathcal{P_D}$, $G \in \mathcal{P_D}$, then $F_{X\vee Y} \in \mathcal{P_D}$.
\item
If $F \in \mathcal{OL}$, $G \in \mathcal{OL}$, then $F_{X\vee Y} \in \mathcal{OL}$.
\item
If $F \in \mathcal{OT}$, $G \in \mathcal{OT}$, then $F_{X\vee Y} \in \mathcal{OT}$.
\item
If $F \in \mathcal{L}$, $G \in \mathcal{L}$, then $F_{X\vee Y} \in \mathcal{L}$.
\item
If $F \in \mathcal{T}$, $G \in \mathcal{T}$, then $F_{X\vee Y} \in \mathcal{T}$.
\end{enumerate}
\ethe

\pr~
From \cite[Lem. 2.2]{ko:tang:2008} we find
\beam \label{eq.KPP.2.9}
\bF_{X\vee Y}(x)=\PP[X\vee Y >x] \sim \PP[X > x] + \PP[Y >x]=\bF(x)+\bG(x)\,,
\eeam
as $\xto$. 
\begin{enumerate}
\item
Let consider $v>1$, then from relation \eqref{eq.KPP.2.9} we obtain
\beao
\limsup_{\xto} \dfrac{\bF_{X\vee Y}(v\,x)}{\bF_{X\vee Y}(x)} &=&\limsup_{\xto} \dfrac{\bF(v\,x) +\bG(v\,x) }{\bF(x)+\bG(x)} \\[2mm]
& \leq& \max\left\{\limsup_{\xto} \dfrac{\bF(v\,x)}{\bF(x)}\,,\;\limsup_{\xto} \dfrac{\bG(v\,x) }{\bG(x)} \right\} < 1\,,
\eeao
where in the pre-last step we used the inequality
\beam \label{eq.KPP.2.10}
\min\left\{ \dfrac ab\,,\;\dfrac cd \right\} \leq \dfrac {a+c}{b+d} \leq \max\left\{ \dfrac ab\,,\;\dfrac cd \right\} \,,
\eeam
and in last step we use the relation $F,\,G \in \mathcal{P_D}$. Hence, $F_{X\vee Y} \in \mathcal{P_D}$.
\item
Let consider $t \neq 0$, then from relations \eqref{eq.KPP.2.9}, \eqref{eq.KPP.2.10} and the fact that $F,\,G \in \mathcal{OL}$ we obtain
\beao
\limsup_{\xto} \dfrac{\bF_{X\vee Y}(x-t)}{\bF_{X\vee Y}(x)}&=&\limsup_{\xto} \dfrac{\bF(x-t)+\bG(x-t)}{\bF(x)+\bG(x)}\\[2mm]
&\leq& \max\left\{\limsup_{\xto} \dfrac{\bF(x-t)}{\bF(x)}\,,\;\limsup_{\xto} \dfrac{\bG(x-t) }{\bG(x)} \right\} < \infty\,.
\eeao
Therefore $F_{X\vee Y} \in \mathcal{OL}$.
\item
This follows directly from combination of 1. and 2. parts and the definition of $\mathcal{OT}$.
\item
Let $t>0$ the inequality
\beam \label{eq.KPP.e}
\liminf_{\xto} \dfrac{\bF_{X\vee Y}(x-t)}{\bF_{X\vee Y}(x)} \geq 1\,,
\eeam
holds always. From the other side because of \eqref{eq.KPP.2.6} and that $F,\,G \in \mathcal{L}$ we find
\beam \label{eq.KPP.z}
\liminf_{\xto} \dfrac{\bF_{X\vee Y}(x-t)}{\bF_{X\vee Y}(x)}&=&\liminf_{\xto} \dfrac{\bF(x-t)+\bG(x-t)}{\bF(x)+\bG(x)}\\[2mm] \notag
&\leq&  \max\left\{\limsup_{\xto} \dfrac{\bF(x-t)}{\bF(x)}\,,\;\limsup_{\xto} \dfrac{\bG(x-t) }{\bG(x)} \right\} =1\,,
\eeam
hence from relations \eqref{eq.KPP.e} and \eqref{eq.KPP.z} we get $F_{X\vee Y} \in \mathcal{L} $.
\item
Follows directly from combination of parts 1. and 4..
~\halmos 
\end{enumerate}

Now we examine the closure property with respect to minimum in the same distribution class, but under another dependence structure. The minimum is not so popular,  however in some application of risk management, as for example in marginal expected shortfall, see \cite{li:2022a}, it is helpful to know whether the marginal class remains after the minimum.

The following dependence structure was introduced in \cite{li:2018b} and it was studied in several  paper including \cite{wang:cheng:yan:2022}, \cite{li:2018b} and \cite{ji:wang:yan:cheng:2023}. We give this condition of dependence for non-negative random variables and in case these random variables are not bounded from below, we need another condition, as follows.

\begin{assumption} \label{ass.KPP.2}
For some distributions  $F,\,G$ of the random variables $X,\,Y$ respectively, with support the interval $\bbr_+$, there exists some constant $C>0$, such that
\beam \label{eq.KPP.3}
\PP[X > x\,\;Y>y] \sim C\,\bF(x)\,\bG(y)\,,
\eeam
as $(x,\,y)\to (\infty,\,\infty)$. Then we say that $X$ and $Y$ are strongly asymptotically independent, symbolically $F,\,G \in SAI$.
\end{assumption} 

\bth \label{th.KPP.2}
Let $X,\,Y$ non-negative random variables with distributions $F,\,G$ respectively. Let us assume that $X,\,Y$ satisfy Assumption \ref{ass.KPP.2}. Then we find
\begin{enumerate}
\item
If $F \in \mathcal{P_D}$, $G \in \mathcal{P_D}$, then $F_{X\wedge Y} \in \mathcal{P_D}$.
\item
If $F \in \mathcal{OL}$, $G \in \mathcal{OL}$, then  $F_{X\wedge Y} \in \mathcal{OL}$.
\item
If $F \in \mathcal{OT}$, $G \in \mathcal{OT}$, then $F_{X\wedge Y} \in \mathcal{OT}$.
\item
If $F \in \mathcal{L}$, $G \in \mathcal{L}$, then  $F_{X\wedge Y} \in \mathcal{L}$.
\item
If $F \in \mathcal{T}$, $G \in \mathcal{T}$, then  $F_{X\wedge Y} \in \mathcal{T}$.
\end{enumerate}
\ethe

\pr~ 
\begin{enumerate}
\item
Let $v>1$, thence
\beao
\limsup_{\xto} \dfrac{\bF_{X\wedge Y}(v\,x)}{\bF_{X\wedge Y}(x)} &=&\limsup_{\xto} \dfrac{\PP[X>v\,x\,,\;Y>v\,x] }{\PP[X>x\,,\;Y>x]} =\limsup_{\xto} \dfrac{C\,\bF(v\,x)\,\bG(v\,x)}{C\,\bF(x)\,\bG(x)}  \\[2mm]
&\leq& \limsup_{\xto} \dfrac{\bF(v\,x) }{\bF(x)}\,\limsup_{\xto}\dfrac{\bG(v\,x)}{\bG(x)} < \limsup_{\xto} \dfrac{\bF(v\,x)}{\bF(x)}< 1\,,
\eeao
where in the last two steps, we take into account that $F,\,G \in \mathcal{P_D}$. 
Therefore we conclude $F_{X\wedge Y} \in \mathcal{P_D}$.
\item
Let us assume $t\neq 0$. Then
\beao
\limsup_{\xto} \dfrac{\bF_{X\wedge Y}(x-t)}{\bF_{X\wedge Y}(x)}&=&\limsup_{\xto} \dfrac{C\,\bF(x-t)\,\bG(x-t)}{C\,\bF(x)\,\bG(x)}\\[2mm]
&&\leq \limsup_{\xto} \dfrac{\bF(x-t)\,\bG(x)}{\bF(x)\,\bG(x)}\,\limsup_{\xto}\dfrac{\bG(x-t)}{\bG(x)}  \leq K \limsup_{\xto} \dfrac{\bF(x-t) }{\bF(x)}< \infty\,,
\eeao
where in the last step we used that $F\in \mathcal{OL}$, while the constant $K \in \bbr_+$ follows from $G\in  \mathcal{OL}$. So we find $F_{X\wedge Y} \in \mathcal{OL}$.
\item
It follows directly from 1. and 2. parts and the definition of $\mathcal{OT}$.
\item
For $t>0$ holds \eqref{eq.KPP.e} and for the bound from above we have
\beao
\limsup_{\xto} \dfrac{\bF_{X\wedge Y}(x-t)}{\bF_{X\wedge Y}(x)} &=&\limsup_{\xto} \dfrac{C\,\bF(x-t)\,\bG(x-t)}{C\,\bF(x)\,\bG(x)} \\[2mm]
&& \leq \limsup_{\xto} \dfrac{\bF(x-t)\,\bG(x)}{\bF(x)\,\bG(x)}\,\limsup \dfrac{\bG(x-t)}{\bG(x)} =\limsup_{\xto} \dfrac{\bF(x-t)}{\bF(x)}= 1\,,
\eeao
therefore  this in combination with relation \eqref{eq.KPP.e} give the result.
\item
This follows directly from part 1. and 4. of the theorem.
~\halmos
\end{enumerate}

In the following results we examine the closure properties of classes $\mathcal{P_D}$, $\mathcal{OA}$ and $\mathcal{OT}$ with respect to convolution of independent random variables. In the first two statements we take into account random variables whose support is the whole real axis $\bbr$, while in the other two the support is the positive half-axis $\bbr_+$. Relevant result on non-negative random variables can be found in \cite[Theorem 2.2]{bardoutsos:konstantinides:2011}, under some heavy conditions.

\bth \label{th.KPP.3.3}
Let $X,\,Y$ be independent random variables with distributions $F,\,G$ respectively.
\begin{enumerate}
\item
Let assume that $F$ and  $G$ have  supports the whole real axis $\bbr$.

a) If $F \in \mathcal{P_D}$, $\bG(x)=o[\bF(x)]$ and $G \in \mathcal{OL}$, then $F*G \in \mathcal{P_D}$.

b) If $F \in \mathcal{OT}$, $\bG(x)=o[\bF(x)]$ and $G \in \mathcal{OL}$, then $F*G \in \mathcal{OT}$.

\item
Let assume that $F$ and  $G$ have supports the positive half-axis $\bbr_+$, and there exists their densities $f$ and $g$ respectively, with $f$ to be of bounded increase, $0<\beta_F < \beta_G$, and 
\beao
\liminf_{\xto}x^q\,\bF(x)>0\,,
\eeao 
for some $q \in [\beta_F,\,\beta_G)$.

a) If $F,\,G \in \mathcal{OA}$, then $F*G \in \mathcal{OA}$.

b) If $F,\,G \in \mathcal{OT}$,  then $F*G \in \mathcal{OT}$.

\item
Let $X_1,\,\ldots,\,X_n $, with $n\geq 2$ be independent random variables with distributions $F_1,\,\ldots,\,F_n$ respectively, supports the positive half-axis $\bbr_+$. We assume that $F_1 \in  \mathcal{OA}$ and moreover
\beam \label{eq.KPP.a}
0<\liminf_{\xto} \inf_{1\leq l \leq n} \dfrac{\bF_l(x)}{\bF_1(x)}\leq\limsup_{\xto} \sup_{1\leq l \leq n} \dfrac{\bF_l(x)}{\bF_1(x)} < \infty\,,
\eeam
where the left inequality can be written as
\beam \label{eq.KPP.b}
\limsup_{\xto} \sup_{1\leq l \leq n} \dfrac{\bF_1(x)}{\bF_l(x)} < \infty\,.
\eeam
Then for the distribution of the sum
\beam \label{eq.KPP.ab}
S_n:=\sum_{i=1}^n X_i\,,
\eeam 
holds $F_{S_n} \in  \mathcal{OA} $.

\item
If $F,\,G$ have as support $\bbr_+$ and either $F \in  \mathcal{L} \cap \mathcal{OA}$ and $G \in  \mathcal{OA}$ or alteratively $G \in  \mathcal{L} \cap \mathcal{OA}$ and $F \in  \mathcal{OA}$, then $F*G \in \mathcal{OA}$.

\item
If $F,\,G$ have as support $\bbr_+$ and either $F,\,G \in  \mathcal{L} \cap \mathcal{OA}$ , then $F*G \in \mathcal{L} \cap \mathcal{OA}$.
\end{enumerate}
\ethe

\pr~
\begin{enumerate}

\item

a) For any $v>1$ we obtain
\beam \label{eq.KPP.3.12} \notag
\overline{F*G}(v\,x)=\int_{-\infty}^{\infty} \bF(v\,x-y)\,G(dy) &=&\left( \int_{-\infty}^{vx-x_0}+\int_{vx-x_0}^{v\,x} + \int_{v\,x}^\infty \right)\bF(v\,x-y)\,G(dy) \\[2mm]
&=:& I_1(v\,x) + I_2(v\,x) + I_3(v\,x)\,,
\eeam
where $x_0$ as found in \eqref{eq.KPP.1.2}. Next we see that
\beao
I_1(v\,x) =\int_{-\infty}^{vx-x_0}\bF(v\,x-y)\,G(dy) \leq C\,v^{-q}\,\int_{-\infty}^{vx-x_0}\bF(x-y)\,G(dy) \,,
\eeao
from which we obtain
\beam \label{eq.KPP.3.13} 
I_1(v\,x) \leq  C\,v^{-q}\,\overline{F*G}(x)\,,
\eeam
where the constants $C$ and $q$ are defined in  \eqref{eq.KPP.1.2}. Further
\beam \label{eq.KPP.3.14} 
I_3(v\,x) \leq  \int_{v\,x}^\infty \,G(dy)=\bG(v\,x)\leq \bG(x)=o\left[\overline{F*G}(x)\right]
\eeam
as $\xto$, where the last step follows from the assumption $\bG(x)=o\left[\bF(x)\right]$. Finally for the second integral we find
\beam \label{eq.KPP.3.15} \notag
I_2(v\,x) &=&  \int_{vx-x_0}^{v\,x}  \bF(v\,x-y)\,G(dy)\leq \bF(0)\,\int_{vx-x_0}^{\infty}\,G(dy) =\bF(0)\,\bG(v\,x-x_0) \\[2mm]
&\leq& \Lambda\,\bG(v\,x)\leq \Lambda\,\bG(x)=o\left[\overline{F*G}(x)\right]\,,
\eeam
as $\xto$, where the constant $\Lambda >0$ follows from that $G \in \mathcal{OL}$. Now from relations \eqref{eq.KPP.3.13},  \eqref{eq.KPP.3.14} and  \eqref{eq.KPP.3.15}, in combination with  \eqref{eq.KPP.3.12}, we find
\beam \label{eq.KPP.3.16} 
\dfrac {\overline{F*G}(v\,x)}{\overline{F*G}(x)} \lesssim C\,v^{-q} +o(1)\,,
\eeam
as $\xto$, thence $F*G \in \mathcal{P_D}$.

b) At first we observe that as far $F \in \mathcal{OT} \subsetneq \mathcal{OL}$ and $G \in \mathcal{OL}$, it follows  $F*G \in \mathcal{OL}$, because of the closure property of $\mathcal{OL}$ with respect to convolution, see for example in \cite[Lemma 3.3]{danilenko:markeviciute:siaulys:2017} or \cite[Proposition 3.12]{leipus:siaulys:konstantinides:2023}. Further, by the inclusion $F\in \mathcal{OT} \subsetneq \mathcal{P_D}$ in the first part, we conclude $F*G \in \mathcal{P_D}$. Hence $F*G \in \mathcal{OT}$.

\item

a) From the fact that $F,\,G \in \mathcal{OA} \subseteq \mathcal{OS}$, and using the closure property of $\mathcal{OS}$ with respect to convolution, see \cite{klueppelberg:1990}, we obtain $F*G \in \mathcal{OS}$. Further from $F,\,G \in \mathcal{OA} \subsetneq \mathcal{P_D}$, taking into consideration our assumptions, through \cite[Theorem 2.2]{bardoutsos:konstantinides:2011}, we find  $F*G \in \mathcal{P_D}$. Thus we have $F*G \in \mathcal{OA}$. 

b) Taking again into account \cite[Lemma 3.3]{danilenko:markeviciute:siaulys:2017}, we obtain the closure property of $\mathcal{OL}$ and from \cite[Theorem 2.2]{bardoutsos:konstantinides:2011} we have the closure property of $\mathcal{P_D}$ as well.

\item

At first, from relations \eqref{eq.KPP.a} and \eqref{eq.KPP.b} follows that $\bF_l(x) \asymp \bF_1(x)$, as $\xto$, for any $1\leq l \leq n$. Therefore, since $F_1  \in  \mathcal{OA}  \subsetneq  \mathcal{OS}$, follows $F_l \in  \mathcal{OS}$ for any $1\leq l \leq n$, see \cite[Lemma 3.1]{watanabe:yamamuro:2010} and from closure property of $ \mathcal{OS}$ with respect to convolution we find $F_{S_n} \in  \mathcal{OS}$, see again \cite[Lemma 3.1]{watanabe:yamamuro:2010}. Now from relation \eqref{eq.KPP.a}, through \cite[Lemma 3]{karaseviciene:siaulys:2023} we get that for any $x\in \bbr$ there exists constant $\widehat{C} \geq 1$ such that
\beam \label{eq.KPP.c} 
\overline{F_{S_n}}(x)\leq \widehat{C}^{n-1}\,\bF_1(x)\,,
\eeam
for any $n\geq 2$. Hence, for some $v>1$ and for any $x>x_0$ we obtain
\beam \label{eq.KPP.d} 
\dfrac{\overline{F_{S_n}}(v\,x)}{\overline{F_{S_n}}(x)}\leq \dfrac{\widehat{C}^{n-1}\,\bF_1(v\,x)}{\overline{F_{S_n}}(x)}\leq \widehat{C}^{n-1}\,\dfrac{\bF_1(v\,x)}{\overline{F_{1}}(x)}\leq\widehat{C}^{n-1}\, C\,v^{-q}:=\widetilde{C}\,v^{-q}\,,
\eeam
where $\widetilde{C}:=\widehat{C}^{n-1}\, C$ and in pre-last step we used relation \eqref{eq.KPP.1.2} for $q \in (0,\,\beta_F)$. Therefore, taking into account \eqref{eq.KPP.d} we establish $F_{S_n} \in \mathcal{P_D}$, which in combination with $F_{S_n} \in \mathcal{OS}$, provides $F_{S_n} \in \mathcal{OA}$.

\item

Taking into consideration that $\mathcal{L} \cap \mathcal{OA}\subsetneq \mathcal{OA} \subsetneq \mathcal{OS}$, we obtain from closure property of $\mathcal{OS}$ with respect to convolution of independent variables that $F*G \in \mathcal{OS}$. Hence, for any $v>1$
\beao
\limsup_{\xto} \dfrac {\overline{F*G}(v\,x)}{\overline{F*G}(x)} &=&\limsup_{\xto}\int_0^{\infty-}\dfrac{\bF(v\,x -y)}{\int_0^{\infty-}\bF(x-z)\,G(dz)}\,G(dy)\\[2mm]
&\leq& \limsup_{\xto}\dfrac{\sup_{0<y<\infty} \bF(v\,x -y)}{\bF(x)} \int_0^{\infty-}\dfrac{1}{\int_0^{\infty-}\,G(dz)}\,G(dy)\\[2mm]
&=& \limsup_{\xto}\sup_{0<y<\infty} \dfrac{\bF(v\,x -y)}{\bF(x)} =\limsup_{\xto}\dfrac{\bF(v\,x )}{\bF(x)}<1\,,
\eeao
where the last two steps follow by the properties of class $\mathcal{L}$ and $\mathcal{P_D}$ respectively. Thus $F*G \in \mathcal{P_D}$ and therefore $F*G \in \mathcal{OA}$. Similarly we can support the alternative argument.
\item
It is true because of previous part 4. and the closure property of class $\mathcal{L}$ with respect to convolution of independent random variables, see \cite{embrechts:goldie:1980}.
~\halmos
\end{enumerate}

\bre \label{rem.KPP.3.1}
We can see that class $\mathcal{P_D}$ has closure property with respect to strong equivalence, namely if $F \in \mathcal{P_D}$ and $\bG(x)\sim c\,\bF(x)$, as $\xto$ for some constant $c>0$, then $G \in \mathcal{P_D}$. Indeed, for any $v>1$ holds
\beao
\limsup_{\xto} \dfrac{\bG(v\,x)}{\bG(x)}=\limsup_{\xto} \dfrac{c\,\bF(v\,x)}{c\,\bF(x)}<1\,.
\eeao 
However, class $\mathcal{P_D}$ has NOT closure property with respect to weak-equivalence, namely from $F \in \mathcal{P_D}$ and $G \asymp F$ is not implied $G \in \mathcal{P_D}$. The classes $\mathcal{OL}$ and $\mathcal{OS}$ enjoy closure properties with respect to strong-equivalence and with respect to weak-equivalence, see for example in \cite[Proposition 3.12]{leipus:siaulys:konstantinides:2023} and in \cite[Lemma 3.1]{watanabe:yamamuro:2010}, respectively. Hence, we can conclude that class $\mathcal{OA}$ and $\mathcal{OT}$ have the closure property with respect to strong-equivalence but NOT with respect to weak-equivalence.
\ere

\bre \label{rem.KPP.3.2}
If the random variables $X,\,Y$  are arbitrarily dependent with distributions $F,\,G \in \mathcal{OT}$, then, for any $p \in (0,\,1)$, hold $(p\,F +(1-p)\,G) \in \mathcal{OT}$. Therefore class $\mathcal{OT}$ has closure property with respect to finite mixture. Indeed this follows by \cite[Proposition 3.12(iii)]{leipus:siaulys:konstantinides:2023} for class $\mathcal{OL}$ and by \cite[Theorem 3.2(2)]{konstantinides:passalidis:2024} for class $\mathcal{P_D}$.
\ere

Let see now the stopped sum
\beam \label{eq.KPP.3.17} 
S_N:=\sum_{i=1}^N X_i\,,
\eeam
with $S_0=0$ and $\{X_i\,,\;i\in \bbn\}$ a sequence of random variables with distributions  $\{F_i\,,\;i\in \bbn\}$ respectively, where $N$ is a discrete random variable with support the set $\bbn_0$ and is independent of the sequence $\{X_i\,,\;i\in \bbn\}$. For the concrete case we consider that $N$ has a bounded support $\{0,\,1,\,\ldots,\,k \}$, for some integer $k \in \bbn$. Let denote by $p_n:=\PP[N=n]$ the probability function of $N$, and put $p_0<1$ and
\beao
\sum_{n=0}^k p_n=1\,.
\eeao 
If the distribution of the distribution of the stopped sum $S_N$ is denoted by $F_{S_N}$, then its tail is represented by
\beao
\bF_{S_n}(x):= \PP[S_N > x] = \sum_{n=1}^k p_n\,\PP[S_n > x]\,,
\eeao
with $S_n$ given by \eqref{eq.KPP.ab}. This construction has many practical applications, for example the $X_i$ can represent the claim sizes during a concrete time interval and $N$ can describe the random multitude of claims during this time. For further reference we mention the works \cite{leipus:siaulys:2012}, \cite{karaseviciene:siaulys:2023}.

\bco \label{cor.KPP.3.1}
Let $X_1,\,X_2,\,\ldots $ independent random variables with distributions $F_1,\,F_2,\,\ldots $ respectively. We assume $F_1 \in \mathcal{P_D} $ and $\bF_i(x)=o[\bF(x)]$, as $\xto$, with $F_i \in \mathcal{OL}$ for any integer $i>1$. If $N$ is a discrete random variable with support bounded from above and independent of $X_1,\,X_2,\,\ldots $, then $F_{S_N} \in  \mathcal{P_D}$.
\eco

\pr~
For case $N=1$ we have directly the result. For $N>1$ we can write
\beam  \label{eq.KPP.3.18} 
\bF_{S_N}(x)=p_1\,\bF_1(x) + \cdots +p_m\,\overline{F_1*\cdots * F_m}(x)\,,
\eeam
for some integer $2\leq m \leq k$. Hence, from relations $F_1 \in \mathcal{P_D}$, $F_2 \in \mathcal{OL}$ and $\bF_2(x)=o[\bF_1(x)]$, as $\xto$, follows $F*G \in \mathcal{P_D}$.

Now for $F_1*\cdots*F_m$ we have again that it belongs to class $\mathcal{P_D}$. Indeed, let assume that $F_1*\cdots *F_{m-1} \in \mathcal{P_D}$, then from $F_m \in \mathcal{OL}$ and $\bF_m(x)=O\left[\bF_1(x)\right]$, we obtain
\beao
\dfrac{\bF_m(x)}{\overline{F_1*\cdots * F_{m-1}}(x)} \leq \dfrac{\bF_m(x)}{\bF_1(x)} \longrightarrow 0\,,
\eeao  
as $\xto$. Hence, from the above we find $\bF_m(x)=o\left[ \overline{F_1*\cdots * F_{m-1}}(x)\right]$ and since $F_m \in  \mathcal{OL}$ and  $F_1*\cdots *F_{m-1} \in \mathcal{P_D}$ are true, they imply $F_1*\cdots*F_m \in \mathcal{P_D}$. Therefore, all the convolutions in \eqref{eq.KPP.3.18} belong to $ \mathcal{P_D}$. Thence, for any $v>1$
\beao
\limsup_{\xto} \dfrac{\bF_{S_N} (v\,x)}{\bF_{S_N} (x)}&=&\limsup_{\xto} \dfrac{p_1\,\bF_1(v\,x) + \cdots +p_m\,\overline{F_1*\cdots * F_m}(v\,x)}{p_1\,\bF_1(x) + \cdots +p_m\,\overline{F_1*\cdots * F_m}(x)}\\[2mm]
&\leq&\max\left\{\limsup_{\xto} \dfrac{\bF_1(v\,x)}{\bF_1(x)},\,\ldots,\,\limsup_{\xto} \dfrac{\overline{F_1*\cdots * F_m}(v\,x)}{\overline{F_1*\cdots * F_m}(x)} \right\} < 1\,,
\eeao 
where the last step is due to the fact that the convolutions belong to class $\mathcal{P_D}$.
~\halmos

In following corollary we find a closure property of randomly stopped sum in class $\mathcal{OA}$.
\bco \label{cor.KPP.2.2}
Let $X_1,\,X_2,\,\ldots $ non-negative independent random variables with distributions $F_1,\,F_2,\,\ldots $ respectively. Let assume that $F_1 \in \mathcal{OA}$ and relations \eqref{eq.KPP.a} and \eqref{eq.KPP.b}, for any $l>1$ with $l \in \bbn$. If $N$ is a discrete random variable with support bounded from above and is independent of $X_1,\,X_2,\,\ldots $, then $F_{S_N} \in \mathcal{OA}$.
\eco

\pr~
For $N=1$, we get directly the result. For $N>1$, we use relation \eqref{eq.KPP.3.15} for some $1\leq m \leq k$, where $k$ is the upper bound of $N$. By Theorem \ref{th.KPP.3.3}3. we have that all the terms in right side of \eqref{eq.KPP.3.15} belong to class $\mathcal{OA}$, namely $F_1 \in \mathcal{OA}$, $F_{S_2} \in \mathcal{OA}$, $\dots$, $F_{S_m} \in \mathcal{OA}$. Whence for any $v>1$
\beao
\limsup_{\xto} \dfrac{\overline{F_{S_N}}(v\,x)}{\overline{F_{S_N}}(x)}&=&\limsup_{\xto} \dfrac{p_1\,\bF_{1}(v\,x)+\cdots+p_m\,\overline{F_{S_m}}(v\,x)}{p_1\,\bF_{1}(x)+\cdots+p_m\,\overline{F_{S_m}}(x)}\\[2mm]
&\leq& \max \left\{ \limsup_{\xto}\dfrac{\bF_{1}(v\,x)}{\bF_{1}(x)},\,\ldots,\, \limsup_{\xto}\dfrac{\overline{F_{S_m}}(v\,x)}{\overline{F_{S_m}}(x)}\right\} <1\,,
\eeao
where the last step is due to $\mathcal{OA} \subsetneq \mathcal{P_D}$. From \cite[Theorem 1]{karaseviciene:siaulys:2023} follows $F_{S_N} \in \mathcal{OS}$. So we have $F_{S_N} \in \mathcal{OA}$.
~\halmos

\section{Further convolution properties in class $\mathcal{A}$ }

In this section we follow similar problems with \cite{embrechts:goldie:1980},\cite{leipus:siaulys:2020} , where was studied the closure property of class $\mathcal{S}$ (and respectively $\mathcal{S(\gamma)}$, for $\gamma\geq 0$) with respect to convolution. In fact we do something similar for class $\mathcal{A} $.
In \cite{konstantinides:leipus:siaulys:2023} the same problem studied for the class of strong subexponential distributions.

\ble \label{lem.KPP.2.1}
Let $X,\,Y$ be independent random variables with distributions $F$ and $G$ respectively, with supports the $\bbr_+$.
\begin{enumerate}
	\item
	If $G \in \mathcal{T}$, $F \in \mathcal{A}$ and 
	\beam \label{eq.KPP.e-h} 
	\limsup_{\xto} \dfrac {\bG(x)}{\bF(x)} < \infty\,,
	\eeam 
	then 
	\beam  \label{eq.KPP.h}
	\overline{F*G}(x) \sim \bF(x) + \bG(x)\,,
	\eeam
	as $\xto$, and moreover $F*G \in \mathcal{P_D}$.
	\item
	If $F_{X \vee Y} \in \mathcal{T}$ and  $F*G \in \mathcal{A}$, then holds \eqref{eq.KPP.h}.
\end{enumerate}
\ele

\pr~
\begin{enumerate}

	\item
	Relation  \eqref{eq.KPP.h} follows from \cite[Lemma 2]{embrechts:goldie:1980} and the fact that $\mathcal{T} \subsetneq \mathcal{A} $ and $ \mathcal{A} \subsetneq \mathcal{S}$. Hence, for any $v>1$ we obtain  
	\beam  \label{eq.KPP.th} \notag
	\limsup_{\xto} \dfrac{\overline{F*G}(v\,x)}{\overline{F*G}(x)} &=& \limsup_{\xto} \dfrac{\bF(v\,x)+\bG(v\,x)}{\bF(x)+\bG(x)}\\[2mm]
	&\leq& \max\left\{\limsup_{\xto} \dfrac{\bF(v\,x)}{\bF(x)}\,,\;\limsup_{\xto} \dfrac{\bG(v\,x) }{\bG(x)} \right\} <1\,,
	\eeam
	where the last step follows from $\mathcal{T} \subsetneq \mathcal{P_D} $ and $ \mathcal{A} \subsetneq \mathcal{P_D}$. Therefore $F*G \in \mathcal{P_D}$.
	
	\item
	Follows directly from \cite[Lemma 1]{embrechts:goldie:1980} and from $\mathcal{T} \subsetneq \mathcal{L} $ and $ \mathcal{A} \subsetneq \mathcal{S}$.~\halmos
\end{enumerate}

The next statement gives a closure property with respect to class $ \mathcal{A}$ with respect to convolution.

\bth \label{th.KPP.2.4}
Let  $X,\,Y$ be independent random variables with distributions $F \in \mathcal{T} $ and $G \in \mathcal{A}$ respectively, with supports the $\bbr_+$. If relation \eqref{eq.KPP.e-h} holds, then we obtain the equivalence: $F \in \mathcal{A}$ if and only if $F*G \in \mathcal{A}$.
\ethe

\pr~
\begin{enumerate}
	
	\item[($\Longleftarrow$)]
	Let $F \in \mathcal{T}$, $G \in \mathcal{T}$, $F*G \in \mathcal{A}$ and relation \eqref{eq.KPP.e-h} holds. Then by \cite[Theorem 1]{embrechts:goldie:1980} follows $F \in \mathcal{S} $, and since $F \in \mathcal{T}$ we find $F \in \mathcal{T}\cap \mathcal{S} = \mathcal{A}$.
	
	\item[($\Longrightarrow$)]
	Let $F \in \mathcal{A}$ and  $G \in \mathcal{A}$ and \eqref{eq.KPP.e-h} be true. Then by \cite[Theorem 1]{embrechts:goldie:1980} we find $F*G \in \mathcal{S}$. Now we apply Lemma \ref{lem.KPP.2.1} to prove $F*G \in \mathcal{P_D}$, that finally gives $F*G \in \mathcal{A}$.~\halmos
\end{enumerate}

\bre \label{rem.KPP.3.3}
In the first part of the proof ($\Longleftarrow$), we use the condition $G \in \mathcal{T}$, which is more general than the condition $G \in \mathcal{A}$ in the second part.
\ere

\bpr \label{pr.KPP.2.1} 
If $F*G \in \mathcal{A}$, $F_{X\vee Y} \in  \mathcal{T}$ and \eqref{eq.KPP.e-h} hold, then $F \in  \mathcal{T}$.
\epr

\pr~
From \cite[Proposition 1]{embrechts:goldie:1980}, due to $\mathcal{A} \subsetneq \mathcal{S}$ and $\mathcal{T} \subsetneq \mathcal{L}$, we find $F \in \mathcal{L}$. Let now assume $F \notin \mathcal{P_D} $, then 
\beao
\lim_{\xto} \dfrac{\bF(v\,x)}{\bF(x)}=1\,,
\eeao
for any $v>1$ and more over
\beam  \label{eq.KPP.2.l} \notag
\limsup_{\xto} \dfrac{\overline{F*G}(v\,x)}{\overline{F*G}(x)} &\leq& \limsup_{\xto} \dfrac{\bF(v\,x)+\bG(v\,x)}{\bF_{X\wedge Y}(x)}=\limsup_{\xto} \dfrac{\bF(v\,x)+\bG(v\,x)}{\bF(x)+\bG(x)}\\[2mm]
&\leq& \max\left\{\limsup_{\xto} \dfrac{\bF(v\,x)}{\bF(x)}\,,\;\limsup_{\xto} \dfrac{\bG(v\,x) }{\bG(x)} \right\} \leq 1\vee B=1\,,
\eeam
with
\beao
B:=\limsup_{\xto} \dfrac{\overline{G}(v\,x)}{\overline{G}(x)} \leq 1\,,
\eeao
where in the second step we used Lemma \ref{lem.KPP.2.1} 2., while in the third step, keeping in mind the independence of $X$ and $Y$ we obtain
\beao
\overline{F*G}(x)=\PP[X\vee Y >x]=\bF(x) + \bG(x) - \bF(x)\,\bG(x) \sim \bF(x) + \bG(x)\,,
\eeao
as $\xto$. Thus, by \eqref{eq.KPP.2.l}, but also the symmetric relation
\beam  \label{eq.KPP.2.k} \notag
\limsup_{\xto} \dfrac{\overline{F*G}(v\,x)}{\overline{F*G}(x)} &\geq& \limsup_{\xto} \dfrac{\bF_{X\vee Y}(v\,x)}{\overline{F*G}(x)}=\limsup_{\xto} \dfrac{\bF(v\,x)+\bG(v\,x)}{\bF(x)+\bG(x)}\\[2mm]
&\geq& \min\left\{\limsup_{\xto} \dfrac{\bF(v\,x)}{\bF(x)}\,,\;\limsup_{\xto} \dfrac{\bG(v\,x) }{\bG(x)} \right\} \geq 1\wedge B\,,
\eeam 
hence, if $G \notin \mathcal{P_D}$ then from \eqref{eq.KPP.2.k} we find
\beao
\limsup_{\xto} \dfrac{\overline{F*G}(v\,x)}{\overline{F*G}(x)} \geq 1\,,
\eeao
and from \eqref{eq.KPP.2.l}
\beao
\limsup_{\xto} \dfrac{\overline{F*G}(v\,x)}{\overline{F*G}(x)} \leq 1\,.
\eeao
But thus we found
\beao
\limsup_{\xto} \dfrac{\overline{F*G}(v\,x)}{\overline{F*G}(x)} = 1\,.
\eeao
that contradicts to assumption $F*G \in \mathcal{A} \subsetneq \mathcal{P_D}$, which implies that  $F \in \mathcal{P_D} $, that means $F \in \mathcal{T}$. 
~\halmos

In the following, we consider that distributions $F$ and $G$ have as support the whole real axis $\bbr$ and the random variables $X$ and $Y$ are independent. The next result
can be found in \cite{leipus:siaulys:2020}.

\bpr[Leipus - {\v S}iaulys] \label{pr.KPP.2.2}
Let $F,\,G \in \mathcal{L}(\gamma)$, with $\gamma \geq 0$. Then the following are equivalent
\begin{enumerate}
\item
	$F*G \in \mathcal{S}(\gamma)$.
\item
	$\overline{F*G}(x) \sim \widehat{G}(\gamma)\,\bF(x) + \widehat{F}(\gamma)\,\bG(x)$, as $\xto$.
\item
	$p\,F + (1-p)\,G \in \mathcal{S}(\gamma)$, for any (or equivalently for some) $p\in (0,\,1)$.
\item
	$F_{X\vee Y} \in \mathcal{S}(\gamma)$.
\end{enumerate}
\epr

\bre \label{rem.KPP.2.4}
We should remark some conclusions. It is well-known that 
\beam \label{eq.KPP.lambd}
\bigcup_{\gamma >0} \mathcal{L}(\gamma) \subsetneq \mathcal{P_D}\,,
\eeam
so it follows 
\beam \label{eq.KPP.mu}
\bigcup_{\gamma >0} \mathcal{S}(\gamma) \subsetneq \mathcal{P_D}\,,
\eeam
see for example \cite[section 2.4]{leipus:siaulys:konstantinides:2023}. Hence, from relations \eqref{eq.KPP.lambd} and \eqref{eq.KPP.mu} we find
\beam \label{eq.KPP.nu}
\left(\bigcup_{\gamma \geq 0} \mathcal{L}(\gamma) \right)\cap \mathcal{P_D} =\mathcal{T}\,, \qquad \left(\bigcup_{\gamma \geq 0} \mathcal{S}(\gamma) \right)\cap \mathcal{P_D} =\mathcal{A}\,,
\eeam
so if instead of condition $F,\,G \in \mathcal{L}(\gamma)$, for any $\gamma \geq 0$ we put $F,\,G \in \mathcal{T}$, the Laplace-Stieltjes transforms $\widehat{F}(\gamma),\,\widehat{G}(\gamma)$ are equal to  unity, and the reason is the inclusion property \eqref{eq.KPP.nu}. Therefore, we study Proposition \ref{pr.KPP.2.1} in case $\gamma = 0$ and there we join the intersection with class $\mathcal{P_D}$.  
\ere

In the next result there exists answer to the following question: What are the necessary and sufficient conditions for the closure property of class $\mathcal{A}$ with respect to convolution of independent random variables.

\bth \label{th.KPP.2.5}
Let $F,\,G \in \mathcal{T}$, then the following are equivalent
\begin{enumerate}
	\item
	$F*G \in \mathcal{A}$.
\item
	$\overline{F*G}(x) \sim \bF(x) + \bG(x)$, as $\xto$.
\item
	$p\,F + (1-p)\,G \in \mathcal{A}$, for any (or equivalently for some) $p\in (0,\,1)$.
\item
	$F_{X\vee Y} \in \mathcal{A}$.
\end{enumerate}
\ethe

\pr~
For the closure property of $\mathcal{S}$, in any equivalence, we employ Proposition \ref{pr.KPP.2.2}, with $\gamma =0$.
\begin{enumerate}
\item[$1. \Rightarrow 2.$] 
	Let assume $F*G \in \mathcal{A}$, then by Proposition  \ref{pr.KPP.2.1} and inclusion $\mathcal{A} \subsetneq \mathcal{S}$, follows $\overline{F*G}(x) \sim \bF(x) + \bG(x)$, as $\xto$.
\item[$2. \Rightarrow 1.$] 
	Let assume $\overline{F*G}(x) \sim \bF(x) + \bG(x)$, as $\xto$. 
	By Proposition \ref{pr.KPP.2.1} we find $F*G \in \mathcal{S}$ and
	\beao
	\limsup_{\xto} \dfrac{\overline{F*G}(v\,x)}{\overline{F*G}(x)} &=& \limsup_{\xto} \dfrac{\bF(v\,x)+\bG(v\,x)}{\bF(x)+\bG(x)}\\[2mm]
	&\leq& \max\left\{\limsup_{\xto} \dfrac{\bF(v\,x)}{\bF(x)}\,,\;\limsup_{\xto} \dfrac{\bG(v\,x) }{\bG(x)} \right\} <1\,,
	\eeao
for any $v>1$, where the last step follows by $F,\,G\in \mathcal{T} \subsetneq \mathcal{P_D}$.
\item[$2. \Rightarrow 3.$]
	Let assume $\overline{F*G}(x) \sim \bF(x) + \bG(x)$, as $\xto$. Then by Proposition \ref{pr.KPP.2.1} follows $p\,F + (1-p)\,G \in \mathcal{S}$ and by \cite[Theorem 3.2]{konstantinides:passalidis:2024}(2) we find that  $\mathcal{P_D}$ has closure property with respect to finite mixture, so $p\,F + (1-p)\,G \in \mathcal{A}$, for any (or equivalently for some) $p\in (0,\,1)$. 
\item[$3. \Rightarrow 2.$]
	Let $p\,F + (1-p)\,G \in \mathcal{A}$, for any (or equivalently for some) $p\in (0,\,1)$. Then, since $ \mathcal{A}\subsetneq \mathcal{S}$, we can apply Proposition \ref{pr.KPP.2.2} to find $\overline{F*G}(x) \sim \bF(x) + \bG(x)$, as $\xto$.
\item[$2. \Rightarrow 4.$]
	Let assume $\overline{F*G}(x) \sim \bF(x) + \bG(x)$, as $\xto$. Then by Proposition \ref{pr.KPP.2.2} follows $F_{X\vee Y} \in \mathcal{S}$ and by Theorem \ref{th.KPP.1} 2., which includes independence also, we find that $\mathcal{P_D}$ has closure property with respect to mamixum, see $\mathcal{T} \subsetneq \mathcal{P_D}$. 
\item[$4. \Rightarrow 2.$]
	Let $F_{X\vee Y} \in \mathcal{A}$. Then, since $ \mathcal{A}\subsetneq \mathcal{S}$, we can apply Proposition \ref{pr.KPP.2.2} to find $\overline{F*G}(x) \sim \bF(x) + \bG(x)$, as $\xto$.~\halmos
\end{enumerate}

Next we present a closure property of randomly stopped sum in \eqref{eq.KPP.3.17}, and this time the counting random variable $N$ is NOT necessarily bounded from above. The corresponding result for the subexponential case can be found in \cite[Theorem 3.37]{foss:korshunov:zachary:2013}, in \cite{denisov:foss:korshunov:2010} and in \cite{karaseviciene:siaulys:2023}. 
\bpr \label{pr.KPP.3.2}
Let $X_1,\,X_2,\,\ldots $ independent and identically distributed random variables with common distribution $F\in  \mathcal{A} $. We assume that random variable $N$ is independent of  $X_1,\,X_2,\,\ldots $, with mean value $\E[N] < \infty$, such that $\E[(1+\delta)^N]< \infty$, for some $\delta >0$. Then
\beam \label{eq.KPP.alpha}
\bF_{S_N}(x) \sim \E[N]\,\bF(x)\,,
\eeam
as $\xto$, and $F_{S_N} \in  \mathcal{A}$
\epr

\pr~
By $\mathcal{A}\subsetneq \mathcal{S}$, relation \eqref{eq.KPP.alpha} is true, see for example \cite{foss:korshunov:zachary:2013} and furthermore  $F_{S_N} \in  \mathcal{S}$, which follows from the closure property of class $\mathcal{S}$ with respect to strong equivalence. Next, we check $F_{S_N} \in  \mathcal{P_D}$. Indeed, for any $v>1$ then
\beao
\limsup_{\xto}\dfrac{\bF_{S_N}(v\,x)}{\bF_{S_N}(x)} = \limsup_{\xto}\dfrac{\E[N]\,\bF(v\,x)}{\E[N]\,\bF(x)}< 1 \,,
\eeao
since $F \in  \mathcal{A} \subsetneq \mathcal{P_D} $. Therefore $F_{S_N} \in  \mathcal{A}$.
~\halmos

\section{Closure properties with respect to convolution product} \label{sec.KP.4}

Here, we examine the closure properties of classes $\mathcal{P_D}$, $\mathcal{OT}$ and $\mathcal{OA}$ with respect to the product convolution. Namely, we consider $X,\,Y$ two independent non-negative random variables, with distributions $F,\,G$ respectively. We wonder if from the assumption $F \in \mathcal{B}$, where $\mathcal{B} \in \{\mathcal{P_D}, \mathcal{OT}, \mathcal{OA}\}$, is implied that the distribution $H$ of the product $X\,Y$ belongs in class $\mathcal{B}$. 

There are several applications in risk theory, in GARCH processes and in infinitely divisible stochastic processes where this question appears. For example we meet such problems in \cite{breiman:1965}, \cite{embrechts:goldie:1980}. However, crucial role on this topic was played by \cite{cline:samorodnitsky:1994}, as it studied this issue in several distributions classes $\mathcal{L}$, $\mathcal{S}$, $\mathcal{D}$ among others. Next, important roles were played by \cite{tang:2006}, \cite{tang:2008}, \cite{xu:cheng:wang:chang:2018}, \cite{cui:wang:2020}. In \cite{tang:2006} we find a sufficient condition for closure property of class $\mathcal{A}$, where, although its not mentioned separately, through the proof of Theorem 2.1 follows, that under the condition
\beam  \label{eq.KPP.3.19} 
\bG(c\,x)=o\left[\bH(x)\right]\,,
\eeam
as $\xto$, for any $c>0$, the class $\mathcal{P_D}$ enjoys the closure property. We consider the non-negative claim sizes, where it was possible to avoid condition \eqref{eq.KPP.3.19}. It worth to mention that recently several attempts to study this problem under dependent random variables, as for example in \cite{cadena:omey:vesilo:2022}, \cite{konstantinides:passalidis:2024}, \cite{yang:wang:leipus:siaulys:2013} among others.
Furthermore, in \cite{konstantinides:leipus:siaulys:2022} we find for random variable $X$, whose support is the whole real axis $\bbr$, that class $\mathcal{OA}$ manifests closure property with respect to convolution product under the condition \eqref{eq.KPP.3.19}, from this aspect, the second part of Theorem \ref{th.KPP.3} extends this result in non-negative case.

\ble \label{lem.KPP.1}
Let $X,\,Y$ non-negative, non-degenerated to zero random variables with distributions $F,\,G$ respectively. 
\begin{enumerate}
\item
Let denote by $H_\vep$ the distribution of the product $X\,(Y\vee \vep)$, and consider some $\delta>0$. If $H_\vep \in \mathcal{P_D}$, for any $\vep \in (0,\,\delta)$, then $H \in \mathcal{P_D}$. 
\item
Let denote by $H_{\vep'}$ the distribution of the product $X\,(Y\vee \vep')$ and consider some $\delta>0$. If $H_{\vep'} \in \mathcal{P_D}$, for any $\vep' \in (\delta,\,\infty)$, then $H \in \mathcal{P_D}$.  
\end{enumerate}
\ele

\pr~
\begin{enumerate}
\item
By \cite[eq. (2.1)]{cline:samorodnitsky:1994} we find that for any $\vep >0$ hold the two-side inequalities
\beao
\PP[Y>\vep]\,\bH_\vep(x) \leq \bH(x) \leq \bH_\vep(x)\,,
\eeao
hence, for any $v>1$ it is true that
\beao
\limsup_{\xto}\dfrac{\bH(v\,x)}{\bH(x)} \leq \limsup_{\xto}\dfrac{\bH_\vep(v\,x)}{\bH_\vep(x)}\,\dfrac 1{\PP[Y>\vep]}  < \dfrac 1{\PP[Y>\vep]} \,,
\eeao
where the prelast step is due to $H_\vep \in \mathcal{P_D}$, thence letting $\vep \downarrow 0$ we have the desired result.
\item
It is easy to see that $\bH_{\vep'}(x) \leq \bH(x)$. Furthermore
\beao
\bH(x) &=& \PP[X\,Y >x\,,\;Y\leq \vep'] + \PP[X\,Y > x\,,\;Y> \vep']\\[2mm]
& \leq& \PP[X\,(Y\wedge \vep') > x\,,\; Y \leq \vep'] + \PP[Y>\vep'] \leq \bH_{\vep'}(x) + \PP[Y>\vep']\,,
\eeao
hence we obtain
\beam \label{eq.KPP.5}
\bH_{\vep'}(x) \leq \bH(x) \leq \bH_{\vep'}(x) + \PP[Y> \vep']\,,
\eeam
which implies that for any $v>1$ holds
\beao
\limsup_{\xto}\dfrac{\bH(v\,x)}{\bH(x)} \leq \limsup_{\xto}\dfrac{\bH_\vep(v\,x)+\PP[Y> \vep']}{\bH_{\vep'}(x)}\,\longrightarrow_{\vep \uparrow 0} \limsup_{\xto}\dfrac {\bH_{\vep'}(v\,x)}{\bH_{\vep'}(x)} <1\,,
\eeao
where in the last step we used $H_{\vep'} \in \mathcal{P_D}$.
\end{enumerate}
~\halmos

\bth \label{th.KPP.3}
Let $X,\,Y$ non-negative, independent random variables with distributions $F,\,G$ respectively. Then
\begin{enumerate}
\item
If $F \in \mathcal{P_D}$, then $H \in \mathcal{P_D}$. 
\item
If $F \in \mathcal{OA}$, then $H \in \mathcal{OA}$.
\item
If $F \in \mathcal{OT}$, then $H \in \mathcal{OT}$.
\end{enumerate}
\ethe

\pr~ 
\begin{enumerate}
\item
Let $H_{\vep,\,\vep'}$ be the distribution of the product $X\,Y$, under the condition that it is greater that zero and less that infinity, namely of the product $X\,(Y\wedge \vep',\,Y\vee \vep)$, with $\vep,\,\vep'$ as introduced in Lemma \ref{lem.KPP.1}. 

Hence, for any $v>1$ we obtain
\beao
\bH_{\vep,\,\vep'}(x)=\int_{1/v}^{x/x_0}\bF\left( \dfrac xy \right)\,G(dy) \geq \dfrac {v^q}{C}\,\int_{1/v}^{x/x_0}\bF\left( \dfrac {v\,x}y \right)\,G(dy)=\dfrac {v^q}{C}\bH{\vep,\,\vep'}(x)\,,
\eeao
where in the second step we take into consideration relation \eqref{eq.KPP.1.2}. Namely, since $F \in \mathcal{P_D}$ we find equivalently $\beta_F >0$, therefore we have
\beam \label{eq.KPP.6}
\dfrac{\bH_{\vep,\,\vep'}(v\,x)}{\bH_{\vep,\,\vep'}(x) } \leq C\,v^{-q}\,.
\eeam
Hence, from \eqref{eq.KPP.6} follows $H_{\vep,\,\vep'} \in \mathcal{P_D}$, as far \eqref{eq.KPP.1.2} is equivalent to $\beta_F>0$, that further is equivalent to $F \in \mathcal{P_D}$. Now, it remains to apply Lemma \ref{lem.KPP.1},1.-2. to find $H \in \mathcal{P_D}$.
\item
The closure property of class $\mathcal{P_D}$ follows from the previous part 1., while the closure property of class $\mathcal{OS}$ is implied by \cite{mikutavicius:siaulys:2023}.
\item
The closure property of class $\mathcal{P_D}$ again follows from the previous part 1., while the closure property of class $\mathcal{OL}$ is due to \cite[Theorem 1]{cui:wang:2020}.~\halmos
\end{enumerate}

In \cite[Theorem 1.3]{xu:cheng:wang:chang:2018} we can find that for non-negative independent random variables $X,\,Y$ with distributions $F,\,G$ respectively, that the implication: $F\in \mathcal{S}$ then $H \in \mathcal{S}$ is equivalent to conditions: either $D[F]=\emptyset$, or $D[F]\neq \emptyset$ and holds
\beam \label{eq.KPP.4.22}
\bG\left( \dfrac xd \right) - \bG\left(\dfrac {x+1}d \right) =o\left[ \bH(x) \right]\,,
\eeam
as $\xto$, for any $d \in D[F]$, where $D[F]$ is the set of positive discontinuity points of distribution $F$. By Theorem \ref{th.KPP.3} we find as consequence the following.

\bco \label{cor.KPP.4.2}
Let $X,\,Y$ non-negative, independent random variables with distributions $F,\,G$ respectively. The implication: $F\in \mathcal{A}$ then $H \in \mathcal{A}$ is equivalent to conditions: either $D[F]=\emptyset$, or $D[F]\neq \emptyset$ and holds
\eqref{eq.KPP.4.22}, for any $d \in D[F]$.
\eco 

\bco \label{cor.KPP.4.3}
Let $X,\,Y$ independent real random variables with distributions $F,\,G$ respectively. If $G(0-)=0$, $G(0)<1$ and $F\in \mathcal{D} \cap \mathcal{P_D}$, then $H \in  \mathcal{D} \cap \mathcal{P_D}$.
\eco

\pr~
It follows directly from Theorem \ref{th.KPP.3} and \cite[Proposition 5.4 (i)]{leipus:siaulys:konstantinides:2023}.
~\halmos

\section{Closure property with respect to minimum in random vectors} \label{sec.KP.5}

In this section we consider two important classes of random vectors, with interest in their closure properties with respect to minimum. The only established class of multivariate heavy-tailed distributions is that of multivariate regular variation. This class was already used in risk theory and risk management, see for example \cite{konstantinides:li:2016}, \cite{yang:su:2023}, \cite{li:2022a} and \cite{li:2022b}.
In bigger distribution classes, there were several attempts for definitions of multivariate subexponential distributions, see \cite{cline:resnick:1992}, \cite{omey:2006}, \cite{samorodnitsky:sun:2016}. Following the definition in \cite{konstantinides:passalidis:2024}, where the distribution tail of the random vector ${\bf X}=(X_1,\,\ldots,\,X_n)$ is given as
\beam \label{eq.KPP.5.1} 
{\bf \bF_t}(x)=\PP[X_1 > t_1\,x,\,\ldots,\,X_n>t_n\,x]\,,
\eeam
for any ${\bf t}=(t_1,\,\ldots,\,t_n) \in (0,\,\infty]^n\setminus \{(\infty,\,\ldots,\,\infty)\}$. The vector ${\bf t}$ plays the role of flexibility in the rate of convergence in this tail. As we see this approach consider the tail only as excesses of all the components, and not the usual tail, where it is given by
\beao
{\bf \bF}({\bf t}x)=1-\PP[X_1 \leq t_1\,x,\,\ldots,\,X_n \leq t_n\,x]\,.
\eeao
Although our approach seems more restrictive, nevertheless it is more immediate with respect to uni-variate case and the same time more practical with respect to closure properties.

The main reason for studying random vectors in the frame of applied probability, is the dependence among their components. According to \cite{konstantinides:passalidis:2024} the n-dimensional dominatedly varying distribution is given by the following.

\bde \label{def.KPP.5.1}
Let a random vector ${\bf X}$, with marginal distributions $F_1,\,\ldots,\,F_n$. We assume that there exists a distribution $F$, such that $F\in \mathcal{D}$ and $\bF_i(x) \asymp \bF(x)$, as $\xto$, for any $i=1,\,\ldots,\,n$. If 
\beam \label{eq.KPP.5.2}
\limsup_{\xto} \dfrac{{\bf \bF_{bt}}(x)}{{\bf \bF_t}(x)}=\limsup_{\xto} \dfrac{\PP[X_1 > b_1\,t_1\,x,\,\ldots,\,X_n>b_n\,t_n\,x]}{\PP[X_1 > t_1\,x,\,\ldots,\,X_n>t_n\,x]} < \infty\,,
\eeam
for any ${\bf t}=(t_1,\,\ldots,\,t_n) \in (0,\,\infty]^n\setminus \{(\infty,\,\ldots,\,\infty)\}$ and any ${\bf b} \in (0,\,1)^n$, then the random vector ${\bf X}$ follows a multi-variate dominatedly varying distribution, symbolically ${\bf \bF_t} \in \mathcal{D}_n$.
\ede

\bde \label{def.KPP.5.2}
Let a random vector ${\bf X}$, with marginal distributions $F_1,\,\ldots,\,F_n$. We assume that $F_i \in \mathcal{D}$, for any $i=1,\,\ldots,\,n$. If 
\beam \label{eq.KPP.5.3}
\limsup_{\xto} \dfrac{{\bf \bF_{vt}}(x)}{{\bf \bF_t}(x)}=\limsup_{\xto} \dfrac{\PP[X_1 >v_1\,t_1\,x,\,\ldots,\,X_n>v_n\,t_n\,x]}{\PP[X_1 > t_1\,x,\,\ldots,\,X_n>t_n\,x]} < 1\,,
\eeam
for any ${\bf t}=(t_1,\,\ldots,\,t_n) \in (0,\,\infty]^n\setminus \{(\infty,\,\ldots,\,\infty)\}$ and any ${\bf v} > {\bf 1}:=(1,\,\ldots,\,1)$, then the random vector ${\bf X}$ follows a multi-variate positively decreasing distribution, symbolically ${\bf \bF_t} \in \mathcal{P_D}_n$.
\ede

These definitions were used in the study of closure properties with respect to scalar product, to convolution and to  randomly stopped sum of random vectors. It is worth to mention that both distribution classes permit arbitrarily dependent components, and under some conditions there exists asymptotic dependence, see \cite{konstantinides:passalidis:2024}. Namely, class $\mathcal{D}_n$ to be asymptotic independent requires only the there exists $F \in \mathcal{D}$, such that $\bF_i(x) \asymp \bF(x)$, as $\xto$, for any $i=1,\,\ldots,\,n$. Indeed, in this case relation \eqref{eq.KPP.5.2} follows immediately. More over, any non-degenerated linear combination of components of ${\bf \bF_t} \in \mathcal{D}_n$, enjoys the uni-variate property $\mathcal{D}$.

The intersection of the two classes $\mathcal{D}_n\cap \mathcal{P_D}_n$ is defined by $F_1,\,\ldots,\,F_n \in \mathcal{D}\cap \mathcal{P_D}$ and relations \eqref{eq.KPP.5.2} and \eqref{eq.KPP.5.3}.

Let us give the minimum of random vectors. We consider two random vectors ${\bf X}_1=(X_{1\,1},\,\ldots,\,X_{1\,n})$ and ${\bf X}_2=(X_{2\,1},\,\ldots,\,X_{2\,n})$, with distributions ${\bf \bF_{t_1}}$ and ${\bf \bF_{t_2}}$, respectively. The tail distribution for the minimum of these two random vectors is
\beam \label{eq.KPP.5.4}
{\bf \bF_{t,\,{\bf X}_1\wedge {\bf X}_2}}(x):=\PP[{\bf X}_1\wedge {\bf X}_2 >{\bf t}\,x]=\PP[ X_{1\,1}\wedge X_{2\,1} >t_1\,x,\,\ldots,\,X_{1\,n}\wedge X_{2\,n} >t_n\,x]\,,
\eeam
for any ${\bf t}=(t_1,\,\ldots,\,t_n) \in (0,\,\infty]^n\setminus \{(\infty,\,\ldots,\,\infty)\}$.  Hence, any multivariate class $\mathcal{B}_n$, has closure property with respect to minimum, if ${\bf \bF_{t,\,{\bf X}_1\wedge {\bf X}_2}} \in \mathcal{B}_n$.

Now, we can introduce the asymptotic independence. Although the components of the vectors are arbitrarily independent, the following condition is related to the dependence structure between the two vectors.

\begin{assumption} \label{ass.KPP.5.1}
Let ${\bf X}_1$ and ${\bf X}_2$ be two random vectors. We say that  ${\bf X}_1$, ${\bf X}_2$ are asymptotically independent if
\beam \label{eq.KPP.5.5}
\PP[{\bf X}_1 >{\bf t}\,x,\,{\bf X}_2 >{\bf t}\,x] \sim \PP[{\bf X}_1 >{\bf t}\,x]\,\PP[ {\bf X}_2 >{\bf t}\,x]\,,
\eeam
as $\xto$, for any  ${\bf t}=(t_1,\,\ldots,\,t_n) \in (0,\,\infty]^n\setminus \{(\infty,\,\ldots,\,\infty)\}$. 
\end{assumption}

Next, we provide a preliminary result.

\ble \label{lem.KPP.5.1}
Let $X_1,\,X_2$ be two random variables with distributions $F_1,\,F_2$ respectively, that satisfy Assumption \ref{ass.KPP.2}.
\begin{enumerate}
\item
If $F_1,\,F_2 \in \mathcal{D}$, then $F_{X_1 \wedge X_2} \in \mathcal{D}$. 
\item
If $F_1,\,F_2 \in \mathcal{D}\cap \mathcal{P_D}$, then $F_{X_1 \wedge X_2} \in \mathcal{D}\cap \mathcal{P_D}$.
\item
If $F_1,\,F_2 \in \mathcal{D}\cap \mathcal{T}$, then $F_{X_1 \wedge X_2} \in \mathcal{D}\cap \mathcal{T}$.
\end{enumerate}
\ele 

\pr~
\begin{enumerate}
\item
At first we find
\beao
\bF_{X_1 \wedge X_2}(x)=\PP[X_1\wedge X_2 >x]=\PP[ X_1 >x,\,X_2>x] \sim C\,\bF_1(x)\,\bF_2(x)\,,
\eeao 
as $\xto$, with some $C>0$. Hence, 
\beao
\limsup_{\xto} \dfrac{\PP[X_1\wedge X_2 >b\,x]}{\PP[X_1\wedge X_2 >x]}&=&\limsup_{\xto} \dfrac{C\,\bF_1(b\,x)\,\bF_2(b\,x)}{C\,\bF_1(x)\,\bF_2(x)}\\[2mm]
&\leq& \limsup_{\xto} \dfrac{\bF_1(b\,x)\,\bF_2(x)}{\bF_1(x)\,\bF_2(x)}\,\limsup_{\xto} \dfrac{\bF_1(b\,x)}{\bF_2(x)}\leq \Lambda\,\limsup_{\xto} \dfrac{\bF_1(b\,x)}{\bF_1(x)}< \infty\,,
\eeao 
where the last two steps follow from $F_1 \in \mathcal{D}$ and $F_2 \in \mathcal{D}$ respectively. Thus $F_{X_1 \wedge X_2} \in \mathcal{D}$.
\item
Follows directly from part 1. and Theorem 2.2, 1..
\item
Follows directly from part 2. and Theorem 2.2, 4., or  from part 1. and Theorem 2.2, 5..~\halmos
\end{enumerate}

\bre \label{rem.KPP.5.1}
Now we ready to formulate a theorem, where for the distributions ${\bf F_{t_1}},\,{\bf F_{t_2}}$ there exist $F_1,\,F_2 \in \mathcal{D}$, such that
\beam \label{eq.KPP.5.6}
 \bF_{1,i}(x) \asymp \bF_1(x)\,,\qquad  \bF_{2,i}(x) \asymp \bF_2(x)\,,
\eeam
as $\xto$, for any $i=1,\,\ldots,\,n$. Only for class $\mathcal{D}_n$, because of its definition, we need the condition $\bF_1(x) \asymp \bF_2(x)$, as  $\xto$. So in the next result, although class $\mathcal{D}_n$ seems to have better properties than $\mathcal{P_D}_n$ or their intersection $\mathcal{D}_n\cap\mathcal{P_D}_n$, here looks more restricted. Let us note that Assumption \ref{ass.KPP.5.1} implies that the components are asymptotic independent, as shown in Assumption \ref{ass.KPP.2} with $C=1$, that means  $F_1,\,F_2$ are asymptotic independent.
\ere

\bth \label{th.KPP.5.1}
Let ${\bf X}_1,\,{\bf X}_2$ be two random vectors with distributions ${\bf F_{t_1}},\,{\bf F_{t_2}}$ respectively, that satisfy Assumption \ref{ass.KPP.5.1}.
\begin{enumerate}
\item
If ${\bf F_{t_1}},\,{\bf F_{t_2}} \in \mathcal{D}_n$ with $F_1(x) \asymp F_2(x)$, as  $\xto$, then ${\bf F_{t,\,X_1 \wedge X_2}} \in \mathcal{D}_n$. 
\item
If ${\bf F_{t_1}},\,{\bf F_{t_2}} \in \mathcal{P_D}_n$, then ${\bf F_{t,\,X_1 \wedge X_2}} \in \mathcal{P_D}_n$. 
\item
If ${\bf F_{t_1}},\,{\bf F_{t_2}} \in \mathcal{D}_n \cap \mathcal{P_D}_n$, then ${\bf F_{t,\,X_1 \wedge X_2}} \in \mathcal{D}_n \cap \mathcal{P_D}_n$. 
\end{enumerate}
\ethe

\pr~
\begin{enumerate}
\item
Independently of distributions of the random vectors, it is true
\beam  \label{eq.KPP.5.7} \notag
{\bf F_{t,\,X_1 \wedge X_2}}(x)&=& \PP[X_{1,1}\wedge X_{2,1}> t_1\,x\,\ldots,\,X_{1,n}\wedge X_{2,n}> t_n\,x]\\[2mm] \notag
&=& \PP[X_{1,1}> t_1\,x,\,X_{2,1}> t_1\,x,\,X_{1,2}\wedge X_{2,1}> t_2\,x,\ldots,\,X_{1,n}\wedge X_{2,n}> t_n\,x]\\[2mm] \notag
&=& \PP[X_{1,1}> t_1\,x,\,X_{2,1}> t_1\,x,\ldots,\,X_{1,n}> t_n\,x,\, X_{2,n}> t_n]\\[2mm] \notag
&\sim& \PP[X_{1,1}> t_1\,x,\ldots,\,X_{1,n}> t_n\,x]\,\PP[X_{2,1}> t_1\,x,,\ldots,\, X_{2,n}> t_n\,x]\\[2mm]
&=& \PP[{\bf X}_1> {\bf t}\,x]\,\PP[{\bf X}_2> {\bf t}\,x]\,,
\eeam
as $\xto$, where in the pre-last step we used Assumption \ref{ass.KPP.5.1}. Let now ${\bf b} \in (0,\,1)^n$, thence
\beam  \label{eq.KPP.5.8} \notag
\limsup_{\xto} \dfrac{{\bf \bF_{b\,t,\,{\bf X}_1\wedge {\bf X}_2}}(x)}{{\bf \bF_{t,\,{\bf X}_1\wedge {\bf X}_2}}(x)}&=&\limsup_{\xto} \dfrac{\PP[{\bf X}_1\wedge {\bf X}_2>{\bf b\,t}\,x]}{\PP[{\bf X}_1\wedge {\bf X}_2>{\bf t}\,x]}\\[2mm] \notag
&=&\limsup_{\xto} \dfrac{\PP[{\bf X}_1>{\bf b\,t}\,x]\,\PP[{\bf X}_2>{\bf b\,t}\,x]}{\PP[{\bf X}_1>{\bf t}\,x]\,\PP[{\bf X}_2>{\bf t}\,x]}\\[2mm] \notag
&\leq& \limsup_{\xto} \dfrac{\PP[{\bf X}_1>{\bf b\,t}\,x]\,\PP[{\bf X}_2>{\bf t}\,x]}{\PP[{\bf X}_1>{\bf t}\,x]\,\PP[{\bf X}_2>{\bf t}\,x]}\,\limsup_{\xto} \dfrac{\PP[{\bf X}_2>{\bf b\,t}\,x]}{\PP[{\bf X}_2>{\bf t}\,x]}\\[2mm] 
&\leq& C\, \limsup_{\xto} \dfrac{\PP[{\bf X}_1>{\bf b\,t}\,x]}{\PP[{\bf X}_1>{\bf t}\,x]}<\infty \,,
\eeam
where in the second step we take into account relation \eqref{eq.KPP.5.7}, while constant $C<\infty$ and the last step follow from ${\bf F_{t_2}} \in \mathcal{D}_n$ and ${\bf F_{t_1}} \in \mathcal{D}_n$ respectively. Therefore relation \eqref{eq.KPP.5.8} correspond to relation \eqref{eq.KPP.5.2} for ${\bf X}_1\wedge {\bf X}_2$ instead of ${\bf X}$.

Now we have to show that there exists distribution $F \in \mathcal{D}$, which satisfies asymptotic $\bF(x) \asymp \bF_j(x)$, as $\xto$, for $j=1,\,2$. As far we are NOT sure whether all the minimums are components of ${\bf X}_1$ or of ${\bf X}_2$, we need the condition $\bF_1(x) \asymp \bF_2(x)$, as $\xto$, that provides the distribution of the minimum $F_{1\wedge2}$, to be weak equivalent with the marginal distributions of the new vector, containing the minimums of the components.
Therefore, since $\bF_1,\, \bF_2$ stem from random variables, which ara asymptotic independent and also $F_{1},\,F_{2} \in \mathcal{D}$, taking into consideration Lemma \ref{lem.KPP.5.1}, 1., in special case of Assumption \ref{ass.KPP.2} for $C=1$, we find $F_{1\wedge2} \in \mathcal{D}$.

\item
We proceed now to the closure property of $\mathcal{P_D}_n$. We begin with relation \eqref{eq.KPP.5.7}, whence
\beam  \label{eq.KPP.5.9} \notag
\limsup_{\xto} \dfrac{{\bf \bF_{v\,t,\,{\bf X}_1\wedge {\bf X}_2}}(x)}{{\bf \bF_{t,\,{\bf X}_1\wedge {\bf X}_2}}(x)}&=&\limsup_{\xto} \dfrac{\PP[{\bf X}_1\wedge {\bf X}_2>{\bf v\,t}\,x]}{\PP[{\bf X}_1\wedge {\bf X}_2>{\bf t}\,x]}\\[2mm] \notag
&=&\limsup_{\xto} \dfrac{\PP[{\bf X}_1>{\bf v\,t}\,x]\,\PP[{\bf X}_2>{\bf v\,t}\,x]}{\PP[{\bf X}_1>{\bf t}\,x]\,\PP[{\bf X}_2>{\bf t}\,x]}\\[2mm] \notag
&\leq& \limsup_{\xto} \dfrac{\PP[{\bf X}_1>{\bf v\,t}\,x]\,\PP[{\bf X}_2>{\bf t}\,x]}{\PP[{\bf X}_1>{\bf t}\,x]\,\PP[{\bf X}_2>{\bf t}\,x]}\,\limsup_{\xto} \dfrac{\PP[{\bf X}_2>{\bf v\,t}\,x]}{\PP[{\bf X}_2>{\bf t}\,x]}\\[2mm] 
&<& \limsup_{\xto} \dfrac{\PP[{\bf X}_1>{\bf v\,t}\,x]}{\PP[{\bf X}_1>{\bf t}\,x]}<1 \,,
\eeam 
for some ${\bf v} >{\bf 1}$, where in last two steps we use ${\bf F_{t_2}} \in \mathcal{P_D}_n$ and ${\bf F_{t_1}} \in \mathcal{P_D}_n$ respectively.

Thence, relation  \eqref{eq.KPP.5.9} provides the \eqref{eq.KPP.5.3} with ${\bf X}_1\wedge {\bf X}_2$ instead of ${\bf X}$. Whence, if $F_{(1\wedge 2),i}$ denotes the distribution of minimum of random variables, following distributions $F_{1,i}$ and $F_{2,i}$, for any $i =1,\,\ldots,\,n$, because of asymptotic independence of $F_{1,i}$ and $F_{2,i}$, through Theorem \ref{th.KPP.2}, 1., in special case of Assumption \ref{ass.KPP.2} for $C=1$, then we find $F_{1\wedge 2,i} \in  \mathcal{P_D}$, for any  $i =1,\,\ldots,\,n$. Thus we conclude ${\bf F_{t,\,{\bf X}_1\wedge {\bf X}_2}} \in \mathcal{P_D}_n$.
 
\item
From Lemma \ref{lem.KPP.5.1}, 2., in special case of Assumption \ref{ass.KPP.2} for $C=1$, we get $F_{1\wedge 2,i} \in  \mathcal{D}\cap \mathcal{P_D}$, for any  $i =1,\,\ldots,\,n$. This in combination of \eqref{eq.KPP.5.8} and \eqref{eq.KPP.5.9} imply ${\bf F_{t,\,{\bf X}_1\wedge {\bf X}_2}} \in \mathcal{D}_n \cap \mathcal{P_D}_n$.~\halmos
\end{enumerate}

\end{document}